\documentclass[12pt]{article}

\usepackage[margin=1in]{geometry}
\usepackage[slantedGreek]{mathpazo}
\usepackage[final]{microtype}
\usepackage{xcolor}

\usepackage{amsmath,amssymb,amsthm,mathtools}
\usepackage{bm}

\usepackage{booktabs}
\usepackage{array}

\usepackage{natbib}
\usepackage[colorlinks=true,linkcolor=blue!50!black,
            citecolor=blue!50!black,urlcolor=blue!50!black]{hyperref}

\theoremstyle{plain}
\newtheorem{theorem}{Theorem}

\theoremstyle{definition}

\theoremstyle{remark}
\newtheorem{remark}[theorem]{Remark}

\DeclareMathOperator*{\E}{\mathbb{E}}
\DeclareMathOperator*{\Prob}{\mathbb{P}}
\DeclareMathOperator{\KL}{KL}

\newcommand{\given}{\,|\,}

\title{\large\bfseries E-Values, Bayes Risk,\\[4pt]
  Dual Role of Markov's Inequality}

\author{
  \makebox[.4\linewidth]{Nicholas G.\ Polson\thanks{Professor of
    Econometrics and Statistics.
    Email: \texttt{ngp@chicagobooth.edu}}}\\[2pt]
  \textit{Booth School of Business}\\
  \textit{University of Chicago}
  \and
  \makebox[.4\linewidth]{Daniel Zantedeschi\thanks{Assistant Professor.
    Email: \texttt{danielz@usf.edu}}}\\[2pt]
  \textit{School of Information Systems}\\
  \textit{University of South Florida}
}

\date{Draft: \today}

\begin{document}
\maketitle

\begin{abstract}
Two approaches to hypothesis testing, e-value testing and
Bayes risk minimisation, both invoke Markov's inequality to
control error probabilities.  They differ in \emph{which}
distribution certifies the unit-moment condition: the null
for Type~I error, the alternative for Type~II error.  The
likelihood ratio is not intrinsically an e-value; it acquires
that status only relative to the experiment under which its
expectation is certified.  This note makes the resulting
role-reversal symmetry explicit, traces its asymptotic
sharpening through the information-theoretic arguments of
\citet{BarronClarke1994}, and situates the duality within
the typed evidence calculus of
\citet{PolsonSokolovZantedeschi2026}.
\end{abstract}


\bigskip

\section{The Duality at a Glance}
\label{sec:duality}

Let $B_{10}(x) = p(x\given H_1)/p(x\given H_0)$ denote the
Bayes factor in favour of~$H_1$.  Two facts follow from the
identity $\int p(x\given H_i)\,dx = 1$:
\begin{equation}\label{eq:e-pair}
  \E_{H_0}[B_{10}] = 1
  \qquad\text{and}\qquad
  \E_{H_1}[B_{01}] = 1,
  \quad B_{01} = 1/B_{10}.
\end{equation}

As a bare likelihood ratio, $B_{10}$ is simply a nonnegative
measurable function on the sample space.  It does not carry
evidential force until it is placed inside an experiment that
certifies its moment condition.  Under the experiment
$(\Omega, \mathcal{F}, P_{H_0})$, the identity
$\E_{H_0}[B_{10}] = 1$ certifies $B_{10}$ as an e-value
and licenses Markov's inequality for Type~I control.
Under $(\Omega, \mathcal{F}, P_{H_1})$, the reciprocal
$B_{01}$ satisfies the same unit-moment condition and is
certified as an e-value for Type~II control.  The object is
the same likelihood ratio, up to inversion, but its
evidential role depends on the governing probability law.
Markov's inequality therefore acts not on an abstract
positive variable, but on a variable equipped with a
specific expectation structure that determines which error
probability it controls.

In the terminology of
\citet{PolsonSokolovZantedeschi2026}, this is a statement
about the (L2)$\to$(L3) bridge of the typed evidence
calculus: the likelihood ratio bridges predictive laws~(L2)
to validity certificates~(L3), but the bridge is typed by
the measure under which the crossing is certified.  The
martingale symmetry noted there, that $B_t$ is a
$P_0$-martingale while $1/B_t$ is a $P_1$-martingale, is
the sequential incarnation of
equation~\eqref{eq:e-pair}.  This note isolates that
symmetry and develops its static and asymptotic
consequences.

Applying Markov under each certifying measure gives the
dual bounds
\begin{equation}\label{eq:dual-markov}
  \Prob_{H_0}\!\bigl(B_{10} \ge 1/\alpha\bigr) \le \alpha,
  \qquad
  \Prob_{H_1}\!\bigl(B_{01} \ge 1/\beta\bigr) \le \beta.
\end{equation}
The first is the familiar e-value bound on the Type~I error
\citep{GrunwaldEtAl2024,Shafer2021}.
The second is the corresponding bound on the Type~II error:
the event $\{B_{10} \le \beta\}$, weak evidence for~$H_1$,
has probability at most~$\beta$ under~$H_1$.

\medskip

Table~\ref{tab:duality} summarises the structure.

\begin{table}[ht]
\centering
\small
\begin{tabular}{@{}lccc@{}}
  \toprule
  \textbf{Framework} & \textbf{E-value} &
    \textbf{Certifying measure} & \textbf{Error controlled} \\
  \midrule
  E-value testing &
    $E = B_{10}$, \; $\E_{H_0}[E]\le 1$ &
    $P_{H_0}$ & Type~I $\le\alpha$ \\[4pt]
  Bayes risk &
    $E^{*} = B_{01}$, \; $\E_{H_1}[E^{*}]\le 1$ &
    $P_{H_1}$ & Type~II $\le\beta$ \\
  \bottomrule
\end{tabular}
\caption{Dual applications of Markov's inequality.
  Numerator and denominator exchange roles while the
  certifying measure switches from~$P_{H_0}$
  to~$P_{H_1}$.  Error control is preserved because the
  reciprocal object remains nonnegative with expectation
  at most one under the appropriately governed experiment.}
\label{tab:duality}
\end{table}

\begin{remark}[Bayes risk threshold and the Fubini decomposition]
\label{rem:bayes-risk}
The Bayes risk framework attends to \emph{both} errors
simultaneously by minimising
$r(\delta) = \pi_0\,c_{\textsc{i}}\,\alpha(\delta)
           + \pi_1\,c_{\textsc{ii}}\,\beta(\delta)$.
The optimal threshold is
$t^{*} = \pi_0\,c_{\textsc{i}} /
  (\pi_1\,c_{\textsc{ii}})$,
trading off the two Markov bounds given the prior~$(\pi_0,\pi_1)$
and loss weights~$(c_{\textsc{i}},c_{\textsc{ii}})$.
This is the static specialisation of the Fubini/Bayes-risk
decomposition developed in
\citet{PolsonSokolovZantedeschi2026}, which shows more
generally that under log-loss the Bayes-optimal critical
region is a threshold rule in the likelihood ratio, with the
cutoff determined by priors and losses rather than by
Markov calibration.
\end{remark}

\begin{remark}[Asymmetry of the certification burden]
\label{rem:scope}
The symmetry in Table~\ref{tab:duality} is formal, not
operational.  The two columns impose different specification
burdens.  The null-side guarantee is structurally light: it
requires only a nonnegative variable~$E$ satisfying
$\E_{H_0}[E]\le 1$, with no reference to any alternative.
This is a central attraction of the e-value framework, as
emphasised by \citet{GrunwaldEtAl2024} and
\citet{Shafer2021}.  The alternative-side guarantee is
heavier: it requires an actual distribution~$P_{H_1}$, or at
least a prior over alternatives, to supply the experiment
relative to which the reciprocal ratio is certified.  The
Type~II column activates only once that governing measure is
specified.  The duality is therefore exact at the level of
the Markov mechanism but asymmetric in what must be committed
in advance.  In the language of the typed calculus, this
asymmetry reflects the distinction between validity-layer
certificates (which require only a null) and decision-layer
optimality (which requires both hypotheses and a loss
structure).
\end{remark}

\section{Composite Alternatives and Mixtures}
\label{sec:composite}

The Markov skeleton of Section~\ref{sec:duality} extends
to composite hypotheses and becomes exponentially tight as
the sample size grows.  The argument follows
\citet{BarronClarke1994} and \citet{ClarkeBarron1990}.

\subsection{The marginal likelihood as an e-value}

For a composite alternative $H_1\colon \theta\in\Theta_1$,
no single likelihood $p(x^n\given\theta)$ serves as the
certifying denominator.  A prior $\pi_1(\theta)$ resolves
this by inducing a \emph{mixture experiment}
$P_{\pi_1}(\cdot) = \int P_\theta(\cdot)\,\pi_1(\theta)\,
d\theta$ with marginal likelihood
\[
  m_1(x^n) = \int p(x^n\given\theta)\,\pi_1(\theta)\,d\theta.
\]
The reciprocal Bayes factor
$B_{01}(x^n) = p(x^n\given\theta_0) / m_1(x^n)$
is then certified as an e-value relative to this mixture
experiment:
\[
  \E_{\pi_1}[B_{01}]
  = \int\!\Bigl[\int \frac{p(x^n\given\theta_0)}{m_1(x^n)}\,
      p(x^n\given\theta)\,dx^n\Bigr]\pi_1(\theta)\,d\theta
  = \int p(x^n\given\theta_0)\,dx^n = 1.
\]
The prior therefore plays a precise structural role: it
supplies the governing measure under which the unit-moment
condition holds.  This is the composite-alternative
counterpart of the Bayesian mixture construction in
\citet{PolsonSokolovZantedeschi2026}, where priors over the
null class induce prior-predictive mixtures as the
certifying law for E-processes; here the same device
operates on the alternative side to certify the reciprocal
Bayes factor.  Markov's inequality gives the composite
Type~II bound
\begin{equation}\label{eq:composite-typeII}
  \Prob_{\pi_1}\!\bigl(B_{01} \ge 1/\beta\bigr) \le \beta.
\end{equation}
This bound is exact at the mixture level, not pointwise
in~$\theta$.  The passage from mixture-level validity to
pointwise control requires asymptotic refinement.

\subsection{Barron--Clarke asymptotics: pointwise sharpening}

\citet{BarronClarke1994} sharpen
\eqref{eq:composite-typeII} to a bound at
\emph{fixed}~$\theta_1$ as $n\to\infty$.  Under regularity
conditions,
\begin{equation}\label{eq:bc-asymp}
  \E_{\theta_1}\!\biggl[
    \log\frac{p(X^n\given\theta_0)}{m_1(X^n)}
  \biggr]
  \;\approx\;
  -n\cdot\KL(\theta_1\|\theta_0)
  + \tfrac{d}{2}\log n
  + \log\frac{\sqrt{\det I(\theta_1)}}{\pi_1(\theta_1)}
  + o(1),
\end{equation}
where $I(\theta_1)$ is the Fisher information matrix and $d$
is the parameter dimension.  The leading term drives
$B_{01}$ to zero exponentially in~$n$ at the KL rate; the
$\tfrac{d}{2}\log n$ correction is the minimax redundancy,
the price of not knowing~$\theta_1$.  Applying Markov to
$B_{01} = e^{-\log B_{10}}$ then yields an exponentially
decaying Type~II bound for each fixed~$\theta_1$.

In the language of Section~\ref{sec:duality}, the
reciprocal Bayes factor is certified at the mixture level
for finite~$n$, while at fixed~$\theta_1$ its typical
magnitude is driven exponentially toward zero at rate
$\KL(\theta_1\|\theta_0)$, up to the redundancy correction.
The KL rate here is the same divergence that governs
evidence growth on the null side:
\citet{PolsonSokolovZantedeschi2026} show that
$(1/n)\log B_{10} \to \KL(P_1\|P_0)$ almost surely
under~$P_1$.  The duality thus extends to rates: KL
divergence controls both how fast evidence accumulates
\emph{for}~$H_1$ and how fast the reciprocal evidential
object decays \emph{against}~$H_0$.  This mirrors the
Sanov/inverse-Sanov duality developed in
\citet{PolsonSokolovZantedeschi2026}, where the same rate
function appears with reversed arguments depending on
whether one measures rarity of empirical distributions or
growth of likelihood-ratio evidence.

\subsection{Three levels of the same skeleton}

The three settings form a nested hierarchy, each
refining the same Markov mechanism with additional
structure.  At the base, a simple alternative provides exact
unit-moment certification under~$P_{H_1}$.  Introducing a
composite alternative replaces~$P_{H_1}$ with the mixture
experiment~$P_{\pi_1}$; the certification remains exact,
but now at the mixture level rather than pointwise.
The Barron--Clarke asymptotics add regularity conditions and
recover pointwise control, showing that the certified object
decays exponentially at each fixed~$\theta_1$.
Table~\ref{tab:three-levels} records this progression.

\begin{table}[ht]
\centering
\small
\renewcommand{\arraystretch}{1.3}
\begin{tabular}{@{}lccl@{}}
  \toprule
  \textbf{Setting} &
  \textbf{Certifying measure} &
  \textbf{Type~II bound} &
  \textbf{Asymptotic rate} \\
  \midrule
  Simple $H_1$ &
    $P_{H_1}$ (exact) & $\beta$ & \textemdash \\
  Composite $H_1$, finite $n$ &
    $P_{\pi_1}$ (exact, mixture) & $\beta$ & \textemdash \\
  Composite $H_1$, $n\to\infty$ &
    $P_{\theta_1}$ (pointwise) &
    exponential decay & $\KL(\theta_1\|\theta_0)$ \\
  \bottomrule
\end{tabular}
\caption{Nested refinement of the Type~II Markov bound.
  Each level instantiates the same Markov skeleton with a
  different certifying measure and degree of sharpness.}
\label{tab:three-levels}
\end{table}

\subsection{Information-theoretic interpretation}

The deeper content of \eqref{eq:bc-asymp} is a statement
about coding regret.  The redundancy
$\log\bigl(p(x^n\given\theta)/m_1(x^n)\bigr)$
measures how many extra bits the Bayes predictor spends
relative to the oracle who knows~$\theta$.
\citet{BarronClarke1994} show that its expectation is
bounded by $\tfrac{d}{2}\log n$, the minimax redundancy.
The Type~II Markov bound is the statement that the log
Bayes factor against the true parameter is unlikely to be
large: in coding terms, the Bayes predictive code is
unlikely to assign much shorter code lengths to the wrong
model than to the correct one.  In the asymptotic regime,
this translates into geometric separation.  The reciprocal
evidential object, certified at the mixture level for
finite~$n$, concentrates exponentially around values that
correctly discriminate against the null, with the KL
divergence governing the rate of that concentration.

\begin{remark}[Dawid's supermartingale and pathwise optimality]
\label{rem:dawid}
The Barron--Clarke redundancy bound holds in expectation.
\citet{Dawid1984,Dawid1987} established a strictly stronger
result using a martingale argument in the prequential
framework.  For any competitor distribution~$Q$ and the
Bayesian mixture
$R(\cdot) = \int P_\theta(\cdot)\,\pi(\theta)\,d\theta$
with $\pi > 0$, the log likelihood ratio
$\log\bigl[Q(x^n)/R(x^n)\bigr]$ is a supermartingale
under~$P_\theta$ for almost every~$\theta$.  By the
martingale convergence theorem, this ratio is bounded
almost surely: no single code~$Q$ can beat the Bayesian
mixture by an unbounded amount on
$P_\theta$-typical sequences.

In the notation of this note, the reciprocal Bayes factor
$B_{01} = p(x^n\given\theta_0)/m_1(x^n)$ is exactly
the likelihood ratio of the null code against the
Bayesian mixture.  Dawid's result therefore gives the
pathwise complement to the finite-sample Markov bound:
not only is $\Prob_{\pi_1}(B_{01} \ge 1/\beta) \le \beta$
for each~$n$, but $B_{01} \to 0$ almost surely
under~$P_\theta$ for almost all~$\theta\in\Theta_1$.
The Markov bound controls the probability of a large
deviation at finite horizon; the supermartingale convergence
theorem ensures that, on the path, the reciprocal
evidential object eventually collapses.
This strengthens the coding-regret interpretation:
the Bayesian mixture is not merely optimal in expectation
\citep{BarronClarke1994,ClarkeBarron1990}, but almost
surely optimal as a universal code
\citep{Dawid1984,Dawid1987}.
\end{remark}

\section{Closing Remarks}
\label{sec:discussion}

The standard presentation of e-values makes Markov's
inequality look like a null-calibration device.  It is not.
It is an error-control mechanism for any nonnegative
statistic once the relevant certifying measure has been
fixed.  E-testing exploits the $H_0$-side of this fact;
Bayes risk exploits the $H_1$-side.  The novelty here is
not a new inequality, but a new reading of an old one:
Markov's inequality is the common operational engine of
both frameworks, and the apparent difference between them
comes entirely from which experiment certifies the
unit-moment condition.  The Barron--Clarke redundancy bound
shows that the alternative-side guarantee tightens at the
KL rate as data accumulate.

This note isolates a single structural observation from the
broader typed evidence calculus of
\citet{PolsonSokolovZantedeschi2026}.  That paper separates
sequential evidence into layers of representation, validity,
and decision, and identifies the likelihood ratio as the
canonical evidence representation under log-loss.  The
present note focuses on one consequence: the same
likelihood ratio serves as a validity certificate on
\emph{both} sides of the testing problem, with the
certifying measure determining which error it controls.
The Barron--Clarke asymptotics then show that the
alternative-side certificate is not merely valid but
asymptotically efficient, with its rate governed by the
same KL divergence that controls evidence growth under
the alternative.  The distinction between validity-layer
polynomial bounds (Markov/Ville at scale~$1/b$) and
efficiency-layer exponential rates (KL at
scale~$e^{-n\cdot\KL}$), emphasised in
\citet{PolsonSokolovZantedeschi2026}, is visible here in
concentrated form: the finite-sample Markov bound on
Type~II error is a validity guarantee, while the
Barron--Clarke sharpening reveals the efficiency geometry
beneath it.

The duality also has an almost-sure face.  Under~$P_1$,
the log-likelihood ratio satisfies
$(1/n)\log B_{10} \to \KL(P_1\|P_0)$ almost surely,
so $B_{10}\to\infty$ and $B_{01}\to 0$ at exponential rate.
The Markov bound on Type~II error, which holds for
every finite~$n$, is thus backed by an a.s.\ convergence
guarantee: the reciprocal Bayes factor is not merely
unlikely to be large, it converges to zero with
probability one.  As noted in Remark~\ref{rem:dawid}, this
pathwise guarantee follows from Dawid's supermartingale
argument \citep{Dawid1984,Dawid1987}: the log likelihood
ratio of any competitor against the Bayesian mixture cannot
diverge, which makes the mixture almost surely optimal as
a universal code.  In the sequential setting,
\citet{DawidVovk1999} show more generally that prequential
likelihood ratios form supermartingales yielding both
anytime-valid certificates and a.s.\ consistency,
providing the natural sequential envelope for the
finite-sample duality developed here.

Two further extensions merit investigation.  First,
e-values compose multiplicatively under optional
stopping \citep{GrunwaldEtAl2024,Shafer2021}, but the Bayes
risk analogue for sequential designs, where the
reciprocal process $B_{01,t}$ must remain a supermartingale
under~$P_1$, is less developed.
Second, a sharpening of the composite Type~II bound using
reverse information projections
\citep{GrunwaldEtAl2024,LarssonRamdasRuf2025}
may yield tighter non-asymptotic analogues
of~\eqref{eq:bc-asymp}.

The point of this note is not that the algebra is
surprising.  Rather, it is that the same likelihood ratio
acquires distinct inferential identities depending on which
experiment certifies its moment condition.  There is a small
historical aptness in this: Markov's own mathematics taught
probability to attend not only to isolated events, but to
sequences whose meaning depends on the governing law.  In
that sense, the conclusion is quietly Markovian: the
inferential role of the object changes with the measure
under which it is read.

\bigskip
\bibliography{duality}

\begin{thebibliography}{9}
\providecommand{\natexlab}[1]{#1}
\providecommand{\url}[1]{\texttt{#1}}
\expandafter\ifx\csname urlstyle\endcsname\relax
  \providecommand{\doi}[1]{doi: #1}\else
  \providecommand{\doi}{doi: \begingroup \urlstyle{rm}\Url}\fi

\bibitem[Barron and Clarke(1994)]{BarronClarke1994}
Andrew~R. Barron and Bertrand~S. Clarke.
\newblock Jeffrey's prior is asymptotically least favorable under entropy risk.
\newblock \emph{Journal of Statistical Planning and Inference}, 41:\penalty0
  37--60, 1994.

\bibitem[Clarke and Barron(1990)]{ClarkeBarron1990}
Bertrand~S. Clarke and Andrew~R. Barron.
\newblock Information-theoretic asymptotics of {B}ayes methods.
\newblock \emph{IEEE Transactions on Information Theory}, 36\penalty0
  (3):\penalty0 453--471, 1990.

\bibitem[Dawid(1984)]{Dawid1984}
A.~P. Dawid.
\newblock Present position and potential developments: Some personal views.
  {S}tatistical theory. {T}he prequential approach.
\newblock \emph{Journal of the Royal Statistical Society: Series A},
  147\penalty0 (2):\penalty0 278--292, 1984.

\bibitem[Dawid(1987)]{Dawid1987}
A.~P. Dawid.
\newblock Discussion of the papers by {D}r {R}issanen and {P}rofessors
  {W}allace and {F}reeman.
\newblock \emph{Journal of the Royal Statistical Society: Series B},
  49\penalty0 (3):\penalty0 253--254, 1987.

\bibitem[Dawid and Vovk(1999)]{DawidVovk1999}
A.~P. Dawid and Vladimir Vovk.
\newblock Prequential probability: Principles and properties.
\newblock \emph{Bernoulli}, 5\penalty0 (1):\penalty0 125--162, 1999.

\bibitem[Gr{\"u}nwald et~al.(2024)Gr{\"u}nwald, de~Heide, and
  Koolen]{GrunwaldEtAl2024}
Peter Gr{\"u}nwald, Rianne de~Heide, and Wouter~M. Koolen.
\newblock Safe testing.
\newblock \emph{Journal of the Royal Statistical Society: Series B},
  86\penalty0 (5):\penalty0 1091--1128, 2024.

\bibitem[Larsson et~al.(2025)Larsson, Ramdas, and Ruf]{LarssonRamdasRuf2025}
Martin Larsson, Aaditya Ramdas, and Johannes Ruf.
\newblock The numeraire e-variable and reverse information projection.
\newblock \emph{The Annals of Statistics}, 2025.
\newblock \doi{10.1214/24-AOS2487}.

\bibitem[Polson et~al.(2026)Polson, Sokolov, and
  Zantedeschi]{PolsonSokolovZantedeschi2026}
Nicholas~G. Polson, Vadim Sokolov, and Daniel Zantedeschi.
\newblock Bayes, e-values, and testing.
\newblock Submitted to the Journal of Machine Learning Research, 2026.

\bibitem[Shafer(2021)]{Shafer2021}
Glenn Shafer.
\newblock Testing by betting: A strategy for statistical and scientific
  communication.
\newblock \emph{Journal of the Royal Statistical Society: Series A},
  184\penalty0 (2):\penalty0 407--431, 2021.

\end{thebibliography}

\end{document}